\documentclass{article}
\def\be{\begin{equation}}
\def\ee{\end{equation}}
\newcommand{\brak}[1]{\left\{  \begin{array}{lllll} #1 \end{array} \right. }
\newcommand{\ff}[1]{{\mbox{\boldmath $#1$}}}
\def\a{\alpha}
\def\b{\beta}
\def\lam{\lambda}
\def\x{\ff{x}}
\def\pk{\ff{p}}
\def\vem{\ff{\varepsilon}}
\newcommand{\mat}[1]{{\left(\begin{array}{cccccc} #1 \end{array}\right)}}

\begin{document}
\title{Swarm Intelligence Based Algorithms: A Critical Analysis }
\author{Xin-She Yang \\
School of Science and Technology, \\ Middlesex University, London NW4 4BT, United Kingdom.
}

\date{}

\maketitle

\begin{abstract}
Many optimization algorithms have been developed by drawing inspiration from swarm intelligence (SI).
These SI-based algorithms can have some advantages over traditional algorithms. In this paper,
we carry out a critical analysis of these SI-based algorithms by analyzing their ways to mimic
evolutionary operators. We also analyze the ways of achieving exploration and exploitation in
algorithms by using mutation, crossover and selection. In addition, we also look at algorithms
using dynamic systems, self-organization and Markov chain framework. Finally, we provide
some discussions and topics for further research.
\end{abstract}

{\bf Citation Detail:} X. S. Yang, Swarm intelligence based algorithms: a critical analysis, \\
{\it Evolutionary Intelligence}, vol. 7, no. 1, pp. 17-28 (2014).

\section{Introduction}

Swarm intelligence (SI) has received great interests and attention in the literature.
In the communities of optimization, computational intelligence and computer science,
bio-inspired algorithms, especially those swarm-intelligence-based algorithms,
have become very popular. In fact, these nature-inspired metaheuristic algorithms
are now among the most widely used algorithms for optimization and computational intelligence
\cite{Glover,Kennedy,KozielYang,Yang2010book,Gandomi}. SI-based algorithms such as ant and bee algorithms,
particle swarm optimization, cuckoo search
and firefly algorithm can have many advantages over conventional algorithms
\cite{Dorigo,Kennedy,YangFA,YangBA2012,Gandomi2}.

An optimization algorithm is an iterative procedure, starting
from an initial guess. After a certain (sufficiently large) number
of iterations, it may converge towards to a stable solution, ideally
the optimal solution to the problem of interest. This is essentially
a self-organizing system with solutions as states, and the converged
solutions as attractors.  Such an iterative, self-organizing system
can evolve, according to a set of rules or mathematical equations.
As a result,  such a complex system can interact and self-organize into certain
converged states, showing some emergent characteristics of self-organization.
In this sense, the proper design of an efficient optimization algorithm
is equivalent to finding efficient ways to mimic
the evolution of a self-organizing system \cite{Ashby,Keller}.

Alternatively, we can view an algorithm as Markov chains, and the behaviour
of an algorithm is controlled by its solution states and transition moves.
Indeed, different views can help to analyze algorithms from different perspectives.
We can also analyze an algorithm in terms of its key components
such as exploration and exploitation or the ways that generate solutions
using evolutionary operators. In this paper, we will review and discuss swarm
intelligence based algorithms from different perspectives.

Therefore, the paper is organized as follows. Section 2 provides a brief critical
analysis of an optimization algorithm, identifying its key components.
Section 3 provides a detailed review of a few swarm intelligence based algorithms
and analyze their key components. Section 4 highlights the key issues in
performance measures for comparing algorithms, and then Section 5
discusses some trends and observations.
Finally, Section 6 draws conclusions briefly.

\section{Critical Analysis of Optimization Algorithms}

An optimization algorithm can be analyzed from different perspectives.
In this section, we will analyze it as an iterative procedure, a self-organization
system, two conflicting components, and three evolutionary operators.

\subsection{Algorithm as an Iterative Process}

Mathematical speaking, an algorithm $A$ is an iterative process, which aims
to generate a new and better solution $\x^{t+1}$ to a given problem from the current
solution $\x^t$ at iteration or time $t$. For example, the Newton-Raphson method
to find the optimal value of $f(\x)$ is equivalent to finding the critical points
or roots of $f'(\x)=0$ in a $d$-dimensional space \cite{Yang2008Comp,Suli}. That is,
\be \x^{t+1} = \x^t -\frac{f'(\x^t)}{f''(\x^t)}=A(\x^t). \ee
Obviously, the convergence rate may become very slow near the optimal point
where $f'(\x) \rightarrow 0$.  Sometimes, the true convergence rate may not
be as quick as it should be. A simple way to improve the convergence is
to modify the above formula slightly by introducing a parameter $p$ as follows:
\be \x^{t+1} = \x^t -p \frac{f'(\x^t)}{f''(\x^t)}, \quad p=\frac{1}{1-A'(\x_*)}. \ee
Here, $\x_*$ is the optimal solution, or a fixed point of the iterative formula.
It is worth pointing out that the optimal convergence of Newton-Raphson's method
leads to an optimal parameter setting $p$ which depends on the iterative
formula and the optimality $\x_*$ of the objective $f(\x)$ to be optimized.

In general, we can write the above iterative equation as
\be \x^{t+1}=A(\x^t,p). \ee
The above iterative formula is valid for a deterministic method; however, in modern
metaheuristic algorithms, randomization is often used in algorithms, and
in many cases, randomization appears in the form of a set of $m$ random variables $\vem=(\varepsilon_1,...,\varepsilon_m)$
in an algorithm. For example, in simulated annealing, there is one random variable,
while in particle swarm optimization \cite{Kennedy}, there are two random variables.

In addition, there are often a set of $k$ parameters in an algorithm. For example,
in particle swarm optimization, there are 4 parameters (two learning parameters,
one inertia weight, and the population size). In general, we can have a vector
of parameters $\pk=(p_1, ..., p_k)$. Mathematically speaking, we can write an
algorithm with $k$ parameters and $m$ random variables as
\be \x^{t+1}=\ff{A}\Big(\x^t, \ff{p}(t), \ff{\varepsilon}(t)\Big), \label{Equ-Alg} \ee
where $\ff{A}$ is a nonlinear mapping from a given solution (a $d$-dimensional
vector $\x^t$) to a new solution vector $\x^{t+1}$.

It is worth pointing out that the above formula (\ref{Equ-Alg}) is for a trajectory-based,
single agent system. For population-based algorithms with a swarm of $n$ agents or solutions,
we can extend the above iterative formula to
the following formula
\be \mat{\x_1 \\ x_2 \\ \vdots \\ \x_n}^{t+1}=A\Big((\x_1^t, \x_2^t, ..., \x_n^t); (p_1, p_2, ..., p_k);
(\epsilon_1, \epsilon_2, ..., \epsilon_m) \Big)  \mat{\x_1 \\ \x_2 \\ \vdots \\ \x_n}^t, \ee
where $p_1, ..., p_k$ are $k$ algorithm-dependent parameters and $\epsilon_1, ..., \epsilon_m$
are $m$ random variables.

This view of algorithm (\ref{Equ-Alg}) is mainly dynamical or functional. It considers  the
functional (\ref{Equ-Alg}) as a dynamical system, which will evolve towards its
equilibrium or attractor states. In principle, the behavoiur of the system
can be described by the eigenvalues of $A$ and its parameters in terms of linear and/or weakly nonlinear
dynamical system theories. However, this does not provide sufficient insight into the
diversity and complex characteristics. Self-organization may provide better insight
as we can see below.

\subsection{Self-Organization in Nature}

Self-organization occurs in nature in a diverse range of physical, chemical and biological systems.
For example, the formation of Rayleigh-B\'enard convection cells,  Belousov-Zhabotinsky oscillator,
and zebra and tiger skin patterns can all be considered as self-organization. A complex, nonlinear, non-equilibrium system will self-organize
under the right conditions. However, as different systems can behave differently, and the conditions
and route to self-organization may be difficult to define. Simple systems such as a swing pendulum
will self-organize towards an equilibrium state, but such a state is usually simple and uninteresting.
On the other hand, highly nonlinear complex systems such as biological systems can have huge diversity in terms of self-organized states, but it is impossible to describe them in a unified theoretical framework.

From the extensive studies of self-organization, it can be summarized that there are certain conditions
that are commonly in  most systems for self-organization to occur \cite{Ashby,Keller}.
For example, the size of the system should be sufficiently large
with a sufficiently high number of degrees of freedom or possible
states $S$. In addition, the system must allow to evolve for a long time from noise and far from equilibrium
states. Most importantly, a selection mechanism must be in place to ensure self-organization is
possible \cite{Ashby,Keller}. That is, the conditions for self-organization in a complex system are:
\begin{itemize}
\item[$\bullet$] The system size is large with a sufficient number of degrees of freedom or states.

\item[$\bullet$] There is enough diversity in the system such as perturbations, noise, edge of chaos, or far from equilibrium.

\item[$\bullet$] The system is allowed to evolve for a long time.

\item[$\bullet$] A selection  mechanism,  or an unchanging law, acts in the system.

\end{itemize}

In other words, a system with states $S$ will evolve towards the self-organized state $S_*$
driven by a mechanism $\alpha(t)$ with a set of parameters $\alpha$.
That is
\be S \stackrel{\alpha(t)}{\longrightarrow} S_*. \ee

\subsection{Algorithms as Self-Organization}

From the self-organization point of view, algorithm (\ref{Equ-Alg}) is a
self-organization system, starting from many possible states $\x^t$ and trying
to converge to the optimal solution/state $\x_*$, driven by the algorithm $A$
in the following schematic representation:
\be f(\x^t) \stackrel{A}{\longrightarrow} f_{\min}(\x_*). \ee

By comparing the conditions for self-organization and characteristics of algorithms, their similarities
and differences are summarized in Table~\ref{table-alg}.
\begin{table}
\begin{center}
\caption{Similarity between self-organization and an optimization algorithm. \label{table-alg} }
\begin{tabular}{|l|l||l|l|}
\hline
Self-organization & Features & Algorithm & Characteristics \\ \hline
Noise, perturbations & diversity & Randomization & escape local optima \\ \hline
selection mechanism & structure & selection & convergence \\ \hline
re-organization & change of states & evolution & solutions \\ \hline
\end{tabular}
\end{center}
\end{table}

However, there are some significant differences between a self-organizing system and
an algorithm. For self-organization, the avenues to the self-organized states may not be
clear, and time is not an important factor. On the other hand, for an algorithm, the
way that makes an algorithm converge is very important, and the speed of convergence is crucial
so that the minimum computational cost can be achieved to reach truly global optimality.

\subsection{Exploration and Exploitation}

Nature-inspired optimization algorithms can also be analyzed from the ways they explore the search space.
In essence, all algorithms should have two key components: exploitation and exploration,
which are also referred to as intensification and  diversification \cite{Blum,Yang2010book}.

Exploitation uses any information obtained from the problem of interest so as to help to
generate new solutions that are better than existing solutions. However, this process is typically
local, and information (such as gradient) is also local. Therefore, it is for local search.
For example, hill-climbing is a method that uses derivative information to guide the
search procedure. In fact, new steps always try to climb up the local gradient. The advantage
of exploitation is that it usually leads to very high convergence rates, but its disadvantage
is that it can get stuck in a local optimum because the final solution point largely depends on
the starting point.

On the other hand, exploration makes it possible to explore the search space more efficiently,
and it can generate solutions with enough diversity and far from the current solutions.
Therefore, the search is typically on a global scale. The advantage of exploration is that it
is less likely to get stuck in a local mode, and the global optimality can be more accessible.
However, its disadvantages are slow convergence and waste of lot computational efforts because
many new solutions can be far from global optimality.

Therefore, a final balance is required so that an algorithm can achieve good performance.
Too much exploitation and too little exploration means the system may converge more quickly,
but the probability of finding the true global optimality may be low. On the other hand,
too little exploitation and too much exploration can cause the search path wonder around
with very slow convergence. The optimal balance should mean the right amount of
exploration and exploitation, which may lead to the optimal performance of an algorithm.
Therefore, balance is crucially important \cite{Blum}.

However, how to achieve such balance is still an open problem. In fact, no algorithm
can claim to have achieved such balance in the current literature. In essence,
the balance itself is a hyper-optimization problem, because it is the optimization
of an optimization algorithm. In addition, such balance may depend on many factors
such as the working mechanism of an algorithm, its setting of parameters, tuning and
control of these parameters and even the problem to be considered. Furthermore,
such balance may not universally exist, and it may vary from problem to problem.
This is consistent with the so-called ``no-free-lunch (NFL)'' theorems \cite{Wolpert}.

The NFL theorems indicated that for any two algorithms $A$ and $B$, if $A$
performs better than $B$ for some problems, there must be some problems on which that
$B$ will perform better than $A$. That is, if measured over all possible problems,
the average performance of both algorithms are essentially equivalent.
In other words, there is no universally better algorithm that can be efficient
for all problems. Though theoretically solid, NFL theorem may have limited impact
in practice because we do not need to solve all problems and we do not need the
average performance either. One of the main aims of problem-solving in practice
is to try to find the optimal or high-quality feasible solution in a  short,
acceptable timescale. For a given type of problem, some algorithms can indeed
perform better than others. For example, for convex problems, algorithms that
use problem-specific convexity information will perform better than black-box
type algorithms. However, some studies suggest that free lunches can exist
for some types of problems, especially for co-evolutionary approaches \cite{Wolpert2}.

These unresolved problems and mystery can motivate more research in this area,
and it can be expected relevant literature will increase in the near future.

\subsection{Evolutionary Operators}

By directly looking at the operations  of an algorithm,
it is also helpful to see how it works. Let us take genetic algorithms as an example. Genetic algorithms (GA) are a class of
algorithms based on the abstraction of Darwin's evolution
of biological systems, pioneered by J. Holland and his collaborators in the 1960s
and 1970s. Genetic algorithms use genetic operators such as
crossover and recombination, mutation, and selection \cite{Holland}. It has been shown that
genetic algorithms have many advantages over traditional algorithms. Three advantages
are distinct: gradient-free, highly explorative, and parallelism. No gradient/derivative
information is needed in GA, and thus GA can deal with complex, discontinuous problems.
The stochastic nature of crossover and mutation make GA explore the search space
more effectively and the global optimality is more likely to be reached. In addition,
genetic algorithms are population-based with multiple chromosomes, and thus it is possible
to implement in a parallel manner \cite{Yang2008,Yang2010book}.

The three key evolutionary operators in genetic algorithms can be summarized as follows:
\begin{itemize}
\item[$\bullet$] Crossover: the recombination of two parent chromosomes (solutions) by exchanging part
of one chromosome with a corresponding part of another so as to produce offsprings (new solutions).

\item[$\bullet$] Mutation: the change of part of a chromosome (a bit or several bits) to generate new
genetic characteristics. In binary string encoding, mutation can be simply achieved by flipping
between 0 and 1. Mutation can occur at a single site or multiple sites simultaneously.

\item[$\bullet$] Selection: the survival of the fittest, which means the highest quality chromosomes
and/characteristics will stay within the population. This often takes some form of elitism,
and the simplest form is to let the best genes to pass onto the next generations in the population.
\end{itemize}

Mathematically speaking, crossover is a mixing process
in a subspace \cite{Bel,Bel2}. This can be seen by an example. For a problem with $8$ dimensions with a total search
space $\Omega=\Re^8$, the parents solutions are drawn from $\x_1=[a a a a a a b b]$, $\x_2=[a a a a a a a a]$
where $a$ and $b$ can be binary digits or a real value for $i$th component/dimension
of a solution to a problem. Whatever the crossover may be, the offsprings will
be one of the 4 possible combinations: $[a a a a a a b a]$, $[a a a a a a a b]$, $[a a a a a a a a]$
and $[a a a a a a b b]$. In any case, the first 6 variables are always $[a a a a a a]$.
This means that the crossover or recombination will lead to a solution
in a subspace where only the 7th and 8th variables are different. No solution will be in the subspace
that the first 6 variables will be different. In other words,
the crossover operator will only create solutions within a subspace $S=[a a a a a a] \cup \Re^2 \subset \Omega$
in the present case. Therefore, crossover can be a local search operator.
However, crossover can also be global if their parents are significantly different
and the subspace is huge. In fact, the subspace for crossover can be extremely large, and thus crossover in
a large subspace can be practically global. For example, if the parent population
is sampled from the whole search space and if the subspace is almost the entire search space,
the crossover acts a global operator. Therefore, care should be taken when we say
a global operator or a local operator, as their boundaries can be rather vague.

On the other hand, mutation provides a mechanism for potential global exploration. In the above example,
if we mutate one of the solution in the first dimension, it will generate a solution
not in the subspace. For example, for a solution $\x_1=[a a a a a a b b] \in S$, if its first $a$ becomes
$b$, then it generates a new solution $\x_q=[b a a a a a b b] \notin S$. In fact, $\x_q$ can be
very difficult from the existing solutions, and can thus jump out of any previous subspace.
For this reason, mutation operator can be considered as an operator with a global search capability.
However,
it is worth pointing out again that all this depends on the step size and rate of mutation.
If the mutation rate is very low, and the solutions generated are not far away
(either limited by low frequency mutation or few sites mutation), mutation may
form a local random walk. In this case, mutation can be essentially local. Again,
whether mutation is global or local, it depends on its frequency and the diversity
of the new solutions.

Both crossover and mutation will provide the diversity for new solutions. However, crossover provides
good mixing, and its diversity is mainly limited in the subspace. Mutation can provide better diversity,
though far-away solutions may drive the population away from converged/evolved characteristics.

It is worth pointing out that selection is a special operator which has a dual role: to choose the
best solutions in a subspace and to provide a driving force for self-organization or convergence.
Without selection, there is no driving force to choose what is the best for the system, and therefore
selection enables a system to evolve with a goal. With a proper selection mechanism, only the fittest solutions
and desired stages may be allowed to gradually pass on, while unfit solutions in the population will
gradually die out. Selection can be as simple as a high degree of elitism, and only the best is selected.
Obviously, other forms of selection mechanisms such as fitness-proportional crossover can be used.

The role and main functions of three evolutionary operators can be categorized as follows:
\begin{itemize}
\item[$\bullet$] Crossover is mainly for mixing within a subspace (that can be huge).
It will help to make the system converge.

\item[$\bullet$] Mutation provides a main mechanism for escaping local modes,
and it can be generalized as a randomization technique.

\item[$\bullet$] Selection provides a driving force for the system to evolve towards the
desired states. It is essentially an intensive exploitative action.

\end{itemize}

Mutation can also take different forms, and one way is to simply use stochastic
moves or randomization. For example,  the traditional Hooke-Jeeves pattern search (PS) is also a gradient-free
method, which has inspired many new algorithms. The key step in pattern
search is the consecutive increment of one dimension followed by the
increment along the other dimensions. The steps will be tried and shrank
when necessary \cite{Hooke}.

Nowadays, the pattern search idea can be generated by the following equation:
\be x_i=x_i+ \delta x_i=x_i+(x_{\rm new move}-x_i), \quad i=1,2, ..., d.
\ee This can be written as a vector equation \be \x=\x+\delta \x=\x +
(\x_{\rm new move}-\x). \label{Hooke-equ} \ee As we can see here, $\delta
\x$ acts as a mutation operator in $2d$ directions in a $d$-dimensional
space.

Differential evolution (DE) was developed by R. Storn and K. Price by their
seminal papers  in 1996 and 1997 \cite{Storn,StornPrice}.
In fact, modern differential evolution (DE) has strong similarity to
the traditional pattern search mutation operator. In fact, the mutation in DE can be
viewed as the generalized pattern search in any random direction
($\x_p-\x_q$) by
\be \x_i=\x_r+F(\x_p-\x_q), \ee where $F$ is the differential
weight in the range of $[0,2]$. Here, $r,p,q,i$ are four different integers
generated by random permutation.

In addition, DE also has a crossover operator which is controlled by a crossover probability $C_r \in [0,1]$
and the actual crossover can be carried out in two ways: binomial and exponential.
Selection is essentially the same as that used in genetic algorithms. It is to select
the most fittest, and for the minimization problem, the minimum objective value.
Therefore, we have
\be \x_i^{t+1} =\brak{\ff{u}_i^{t+1} & & \textrm{if } \; f(\ff{u}_i^{t+1}) \le f(\x_i^t),  \\[10pt]
\x_i^t & \textrm{otherwise.}} \label{de-equ-555} \ee
Most studies have focused on the choice of $F$, $C_r$, and $n$ as well as
the modification of the mutation scheme. In addition,
it can be clearly seen that selection is also used when the condition in the above
equation is checked. Various variants of DE all use crossover, mutation and selection,
and the main differences are in the steps of mutation and crossover. For example, DE/Rand/1/Bin
use three vectors for mutation and binomial crossover \cite{Price}.

\subsection{Role of the Current Best Solution}

The role of current best solution can be two-fold. One the one hand, it seems that the use of the
current global best $\ff{g}^*$ found so far during the iterative search process can speed up the
convergence as it exploits the information during the
search process; however, both observations and mathematical analysis suggest that
the use of the current best solution may somehow lead to  premature convergence
because the converged current best is not necessarily the true global optimal solution
to the problem under consideration. Therefore, the current best solution
may act as a double-edged sword; it can work well once properly used, but it can
also lead to prematured local solutions.

On the other hand, the current best solution can be considered
as a (narrow) type of elitism, which will increase the local
selection or evolution pressure, and can thus help the system to converge.
However, too much pressure can lead to the inability to adapt to
rapidly changing environments, limiting the diversity of the population.

Therefore, mathematical analysis and practical observations seem to suggest
that the use of the current can speed up convergence under
the right conditions, but it will inhibit the ability and probability
of finding the global optimality. Algorithms that use $\ff{g}^*$ may not
have guaranteed global convergence. Consequently, it is no surprise that
metaheuristic algorithms such as particle swarm optimization
and others do not have guaranteed global convergence. However, algorithms do not
use $\ff{g}^*$ may have global convergence. In fact, simulated annealing and cuckoo search
have guaranteed global convergence.

\section{Analysis of Swarm Intelligence Based Algorithms}

In the rest of this paper, we will first briefly introduce the key steps
of a few SI-based algorithms and then analyze them in terms of evolutionary operators
and the ways for exploration and exploitation.

\subsection{Ant and Bee Algorithms}

Ant algorithms, especially the ant colony optimization developed by M. Dorigo \cite{Dorigo,Dorigo2},
mimic the foraging behavior of social ants.
Primarily, all ant algorithms use pheromone as a chemical messenger
and the pheromone concentration as the indicator of quality solutions to
a problem of interest. From the implementation point of view, solutions are
related to the pheromone concentration, leading to routes and paths marked
by the higher pheromone concentrations as better solutions, which can be suitable for some problems
such as discrete combinatorial problems.

Looking at ant colony optimization closely, random route generation is
primarily mutation, while pheromone-based selection provides a mechanism
for selecting shorter routes. There is no explicit crossover in ant algorithms.
However, mutation is not simple action as flipping digits as in genetic algorithms,
new solutions are essentially generated by fitness-proportional mutation.
For example,  the probability of ants in a network problem
at a particular node $i$ to choose the route from node $i$ to node $j$ is given by
\be p_{ij} =\frac{\phi_{ij}^{\alpha} d_{ij}^{\beta}}
{\sum_{i,j=1}^n \phi_{ij}^{\alpha} d_{ij}^{\beta}}, \label{equ-ant-500} \ee
where $\alpha>0$ and $\beta>0$ are the influence parameters, and $\phi_{ij}$ is the pheromone
concentration on the route between $i$ and $j$. Here, $d_{ij}$ is the desirability
of the same route. The selection is subtly related to some {\it a priori} knowledge about
the route such as the distance $s_{ij}$ that is often used so that $d_{ij} \propto 1/s_{ij}$.

On the other hand, bee algorithms do not usually use pheromone \cite{Nak}.
For example, in the artificial bee colony (ABC) optimization algorithm \cite{Kara},  the bees in a colony
are divided into three groups: employed bees (forager bees),
onlooker bees (observer bees), and scouts.
Randomization is carried out by scout bees and employed bees, and  both are
mainly mutation. Selection is related to the honey or objective.
Again, there is no explicit crossover.

Both ACO and ABC only use mutation and fitness-related selection, they can have
good global search ability. In general, they can explore the search space relatively
effectively, but convergence may be slow because it lacks crossover, and thus
the subspace exploitation ability is very limited. In fact, the lack of crossover
is very common in many metaheuristic algorithms.

In terms of exploration and exploitation, both ant and bee algorithms have strong
exploration ability, but their exploitation ability is comparatively low.
This may explain why they can perform reasonably well for some tough optimization,
but the computational effort such as the number of function evaluations is typically high.

\subsection{Particle Swarm Optimization}

Particle swarm optimization (PSO) was developed by Kennedy and
Eberhart in 1995 \cite{Kennedy}, based on the swarming behaviour such
as fish and bird schooling in nature. In essence, the position and velocity of a particle, $\x_i$ and $\ff{v}_i$, respectively,
can be updated as follows:
\be \ff{v}_i^{t+1}= \ff{v}_i^t  + \a \ff{\epsilon}_1
[\ff{g}^*-\x_i^t] + \b \ff{\epsilon}_2 [\x_i^*-\x_i^t],
\label{pso-speed-100}
\ee
\be \x_i^{t+1}=\x_i^t+\ff{v}_i^{t+1}, \label{pso-speed-140} \ee
where $\ff{\epsilon}_1$ and $\ff{\epsilon}_2$ are two random vectors, and each
entry can take the values between 0 and 1.
The parameters $\a$ and $\b$ are the learning parameters or
acceleration constants, which can typically be taken as, say, $\a\approx \b \approx 2$.

By comparing the above equations with the pattern search in Section 2.3,
we can see that the new position is generated by pattern-search-type
mutation, while selection is implicitly done by using the current global
best solution $\ff{g}^*$ found so far, and also through the individual best
$\x_i^*$. However, the role of individual best is not quite clear, though
the current global best seems very important for selection, as this is shown
in the accelerated particle swarm optimization \cite{Kennedy,Yang2008,YangAPSO}.

Therefore, PSO consists of mainly mutation and selection. There is no crossover
in PSO, which means that PSO can have a high mobility in particles with a high
degree of exploration. However, the use of $\ff{g}^*$ seems strongly selective,
which may be like a double-edge sword. Its advantage is that it helps to speed up
the convergence by drawing towards the current best $\ff{g}^*$, while at the same time
it may lead to premature convergence even though this may not be the true optimal solution
of the problem of interest.

\subsection{Firefly Algorithm}

Firefly Algorithm (FA) was developed by Xin-She Yang in 2008 \cite{Yang2008,YangFA,YangFA2},
which was based on the flashing patterns and behaviour of tropical fireflies.
FA is simple, flexible and easy to implement.

The movement of a firefly $i$ is attracted to another more attractive (brighter) firefly $j$ is determined by
\be    \x_i^{t+1} =\x_i^t + \beta_0 e^{-\gamma r^2_{ij} } (\x_j^t-\x_i^t) + \alpha \; \ff{\epsilon}_i^t, \label{FA-equ} \ee
where the second term is due to the attraction, and $\beta_0$ is the attractiveness at $r=0$.
The third term is randomization with $\alpha$ being the randomization parameter, and $\ff{\epsilon}_i^t$ is a vector of
random numbers drawn from a Gaussian distribution at time $t$. Other studies also use
the randomization in terms of $\ff{\epsilon}_i^t$ that can easily be extended to other distributions such as L\'evy flights \cite{YangFA,YangFA2}.
A comprehensive review of firefly algorithm and its variants has been carried out by Fister et al. \cite{Fister}.

From the above equation, we can see that mutation is used for both local and global search.
When $\ff{\epsilon}_i^t$ is drawn from a Gaussian distribution and L\'evy flights, it produces
mutation on a larger scale. On the other hand, if $\alpha$ is chosen to be a very small value,
then mutation can be very small, and thus limited to a local subspace. Interestingly, there is no explicit
selection in the formula as $\ff{g}^*$ is not used in FA.  However, during the update step in the loops
in FA, ranking as well as selection is used.

One novel feature of FA is that attraction is used, and this is the first of its kind in
any SI-based algorithms. Since local attraction is stronger than long-distance attraction,
the population in FA can automatically subdivide into multiple subgroups, and each group
can potentially swarm around a local mode. Among all the local mode, there is always a global
best solution which is the true optimality of the problem. Thus, FA can deal with multimodal problems naturally
and efficiently.

From Eq.~(\ref{FA-equ}), we can see that FA degenerates into a variant of differential evolution
when $\gamma=0$ and $\alpha=0$. In addition, when $\beta_0=0$, it degenerates into simulated annealing (SA).
Further, when $\x_j^t$ is replaced by $\ff{g}^*$, FA also becomes the accelerated PSO.
Therefore, DE, APSO and SA are special cases of the firefly algorithm, and thus FA can have the
advantages of these algorithms. It is no surprise that FA can be versatile and efficient,
and perform better than other algorithms such as GA and PSO.

As a result, the firefly algorithm and its variants have been applied in a diverse range
of applications \cite{Azad,Fara,Hass,Horng,Jati,Nandy,Palit}, including hard problems and
multiobjective problems \cite{Jati,Sayadi,Senthil,Sriv,Yousif,Zaman}.

\subsection{Cuckoo Search}

Cuckoo search (CS) was developed in 2009 by Xin-She Yang and Suash Deb \cite{YangDeb}.
CS is based on the brood parasitism of some cuckoo species. In addition, this
algorithm is enhanced by the so-called L\'evy flights \cite{Pav},
rather than by simple isotropic random walks.
Recent studies show that CS is potentially far more efficient than PSO
and genetic algorithms \cite{YangDeb2010,YangDeb2013,Walton,Chand,Dhiv,Dhiv2,Durgun}.
In fact, CS and its variants have applied in almost every area of engineering design
and applications \cite{Layeb,Mora,Nogh,Sriv,Valian,Yildiz}.

CS uses a balanced combination of a local random walk and the global
explorative random walk, controlled by
a switching parameter $p_a$. The local random walk can be written as
\be \x_i^{t+1}=\x_i^t +\alpha s \otimes H(p_a-\epsilon) \otimes (\x_j^t-\x_k^t), \label{CS-equ1} \ee
where $\x_j^t$ and $\x_k^t$ are two different solutions selected randomly by random permutation,
$H(u)$ is a Heaviside function, $\epsilon$ is a random number drawn from a uniform distribution, and
$s$ is the step size. Here, $\otimes$ is an entry-wise multiplication.

On the other hand, the global random walk is carried out by using L\'evy flights
\be \x_i^{t+1}=\x_i^t+\a L(s,\lam), \label{CS-equ2} \ee
where \be L(s,\lam)=\frac{\lam \Gamma(\lam) \sin (\pi \lam/2)}{\pi}
\frac{1}{s^{1+\lam}}, \quad (s \gg s_0>0). \ee
Here $\a>0$ is the step size scaling factor, which should be related to the scales of the problem of
interest.

CS has two distinct advantages over other algorithms such as GA and SA,
and these advantages are: efficient random walks and balanced mixing.
Since L\'evy flights are usually far more efficient than any other random-walk-based randomization techniques,
CS can be efficient in global search. In fact, recent studies show that CS can have
guaranteed global convergence \cite{Wang}. In addition, the similarity between eggs
can produce better new solutions, which is essentially fitness-proportional generation with a good mixing
ability. In other words, CS has varying mutation realized by L\'evy flights,
and the fitness-proportional generation of new solutions based on similarity
provides a subtle form of crossover. In addition, selection is carried out
by using $p_a$ where the good solutions are passed onto the next generation, while
not so good solutions are replaced by new solutions. Furthermore, simulations also show
that CS can have auto-zooming ability in the sense that new solutions can automatically zoom into the region where
the promising global optimality is located.

In addition, equation (\ref{CS-equ2}) is essentially simulated annealing in the
framework of Markov chains. In Eq.~(\ref{CS-equ1}), if $p_a=1$ and $\alpha s \in [0,1]$,
CS can degenerate into a variant of differential evolution.  Furthermore, if we replace
$\x_j^t$ by the current best solution $\ff{g}^*$, then (\ref{CS-equ1}) can further degenerate into
accelerated particle swarm optimization (APSO) \cite{YangAPSO}. This means that
SA, DE and APSO are special cases of CS, and that is one of the reasons why CS is so efficient.

It is worth pointing out that both CS and FA can capture the main characteristics of
SA, DE and APSO, but there are some significant difference between FA and CS.
One major difference between FA and CS is that FA uses the distance-based, landscape-modulated
attraction. As local attraction is stronger than long-distance attraction, FA can subdivide
the whole population into such populations, while CS does not. In addition,
FA uses random walks and ranking, while CS uses L\'evy flights and random permutation,
which will result in different behaviour. These differences make both FA and CS unique.

\subsection{Bat Algorithm}

The metaheuristic bat algorithm (BA) was developed by Xin-She Yang in 2010 \cite{YangBA}.
It was inspired by the echolocation behavior of microbats. It is the first algorithm
of its kind to use frequency tuning.  Each bat is associated with a
velocity $\ff{v}_i^t$ and a location $\x_i^t$,
at iteration $t$, in a $d$-dimensional search or solution space. Among
all the bats, there exists a current best solution $\x_*$. Therefore,
the above three rules can be translated into the updating equations
for $\x_i^{t}$ and velocities $\ff{v}_i^{t}$:
\be f_i =f_{\min} + (f_{\max}-f_{\min}) \beta, \label{f-equ-150} \ee
\be \ff{v}_i^{t} = \ff{v}_i^{t-1} +  (\x_i^{t-1} - \x_*) f_i , \ee
\be \x_i^{t}=\x_i^{t-1} + \ff{v}_i^t,  \label{f-equ-250} \ee
where $\beta \in [0,1]$ is a random vector drawn from a uniform distribution.

The loudness and pulse emission rates are regulated by the following equations:
\be A_i^{t+1}=\alpha A_{i}^{t}, \ee
and
\be r_i^{t+1}= r_i^0 [1-\exp(-\gamma t)], \label{rate-equ-50} \ee
where $0<\alpha<1$ and $\gamma>0$ are constants. In essence, here $\alpha$ is similar
to the cooling factor of a cooling schedule in simulated annealing.

BA has been extended to multiobjective bat algorithm (MOBA) by Yang \cite{YangBA2011,YangBA2012},
and preliminary results suggested that it is very efficient.

In BA, frequency tuning essentially acts as mutation, while selection pressure
is relatively constant via the use of the current best solution $\x_*$ found so far.
There is no explicit crossover; however, mutation varies due to the variations of loudness
and pulse emission. In addition, the variations of loudness and pulse emission rates also
provide an autozooming ability so that exploitation becomes intensive as the search
is approaching global optimality.

In principle, we can use the same procedure to analyze the key components
and evolutionary operators in all other algorithms such as the most recent
flower pollination algorithm \cite{Yang2013FPA}. Due to the length limit, we will
not analyze other algorithms here. Instead, we will highlight more important issues
such as performance measures and some open problems.

\section{Performance Measures}

Despite the huge number of studies about various metaheuristic algorithms, it still lacks
a good performance measure for comparing different algorithms. In general, there
are two major ways to compare any two algorithms \cite{Yang2010book}.

One way is to compare the accuracy of two algorithms to solve the same problem
for a fixed number of function evaluations.
In most case, the problem to be solved is usually a test function whose minimum
value is often zero by proper formulation. For example, if algorithm $A$ obtains
$0.001$ for, say, $N=1000$ function evaluations, while algorithm $B$ obtains $0.002$
in a run, one tends to say $A$ obtained a better value than $B$.
However, care should be taken when making a statement. Due to the stochastic
nature, multiple runs should be carried out so that meaningful statistical measures
can be obtained. Suppose, we run each algorithm 100 times, and we can get a mean
function value and its standard deviation. For example, if we run algorithm $A$ for
100 independent runs, we get a mean $\mu_A=0.001$ and a standard deviation $\sigma_A=0.01$.
We can write the results as
\be f_A=\mu_A \pm \sigma_A=0.001 \pm 0.01. \label{A-mean} \ee
Similarly, if we run $B$ 100 times, and we get
\be f_B=\mu_B \pm \sigma_B=0.001 \pm 0.02. \label{B-mean} \ee
As the standard deviation $\sigma_B$ is large number than $\sigma_A$,
and their means are the same, we can still say algorithm $A$ obtain
better results because it has smaller or narrower range of confidence interval, though such
improvement of $A$ over $B$ is marginal.

Another way is to compare the numbers of function evaluations required by
two different algorithms for a given accuracy. Suppose, the accuracy $\delta=10^{-3}$
is fixed, if algorithm $A$ requires about $1000$ function evaluations,
while algorithm $B$ requires $1400$ function evaluations, we may say $A$ is better than $B$.
Again, care should be taken when drawing such conclusions. Multiple runs
are needed to ensure meaningful statistical measures. Suppose, we run each
algorithm 100 times independently, we get a mean $1000$ and a corresponding
standard deviation of $300$. We can write the results as
\be N_A=1000 \pm 300.\ee
Similarly, for algorithm $B$, we may get
\be N_B=1400 \pm 300. \ee
In this case, we may draw the conclusion that $A$ is better than $B$
with some statistical confidence.

A third way is to compare the execution times for two algorithms.
However, this approach is not recommended because it has many drawbacks.
Firstly, one has to ensure each algorithm is implemented correctly and efficiently.
Even two algorithms are equally efficient, but the ways of implementing the
calculations can affect the execution time. Secondly, different computers
may have different configurations, and one has to ensure the computation is
carried out on the same computer. Thirdly, even the same computer is used,
the hidden processes (especially true for the Windows operating system)
can consume a lot of computing power at different times, even the computer
is virus free. Finally, even all the above has ensured that the results are good,
such execution times are usually not repeatable by other users.

A fourth way is to use normalize results for comparison, and this approach
is misleading and can produce wrong conclusions \cite{GGWang}. The idea is to choose
an algorithm as the basis of comparison, and normalize the results by others
to produce a ratio. For example, if we use the results by $B$ as the basis,
we have to normalize $A$ to produce
\be \frac{\mu_A}{\mu_B}=1.00. \ee
We can say both $A$ and $B$ are seemingly equally efficient. However, extreme
care must be taken when doing this because we have to consider uncertainty propagations carefully.

As we know for a function of $A/B$ when $A$ and $B$ has no correlation (i.e., their correlation coefficient is zero),
the overall or combined standard deviation is given by
\be \sigma_{A/B}=\frac{A}{B} \sqrt{\frac{\sigma_A^2}{A^2} +\frac{\sigma_B^2}{B}}. \ee
For the results given in Eqs.(\ref{A-mean}) and (\ref{B-mean}), we have
\be \mu_A=0.001, \;\; \mu_B=0.001, \quad \sigma_A=0.01, \;\; \sigma_B=0.02, \ee
which gives
\be \sigma_{A/B}=\frac{0.001}{0.001} \sqrt{ \frac{0.01^2}{0.001^2}+\frac{0.02^2}{0.001^2}} \approx 22.3607. \ee
Therefore, we have the performance of $A$ relative to $B$
\be \mu_{A/B} \pm \sigma_{A/B}=1.00 \pm 22.36,  \ee
which gives the impression that $A$ is worse than $B$ because of the large standard deviation.
This conclusion contradicts with the earlier results.
Even for the special case when $\mu_A=\mu_B=0.001$ and $\sigma_A=\sigma_B=0.02$, we may
have $1.00 \pm 20 \sqrt{2} =1.00 \pm 28.284$, which gives an impression of an even worse solution,
though in reality both algorithm perform almost equally well.
This clearly demonstrate the danger of this approach. Obviously, results should never be normalized
in this ways for comparison \cite{GGWang}. Otherwise, such misleading comparison can render the results invalid.

As we have seen, statistically robust performance measures are yet to be developed,
and we hope this review can inspire more research in this area.

\section{Discussions}

Many optimization algorithms  are based on swarm intelligence, and use population-based approaches.
Most will use some sort of three key evolutionary operators: crossover, mutation and
selection; however, almost all algorithms use mutation and selection, while crossover may
appear in some subtle way in some algorithms. Crossover is efficiently in exploitation
and can often provide good convergence in a subspace. If this subspace is where
the global optimality lies, then crossover with elitism can almost guarantee to
achieve global optimality. However, if the subspace of crossover is not in the region
where the global optimality lies, there is a danger for premature convergence.

The extensive use of mutation and selection can typically enable a stochastic algorithm to
have a high ability of exploration.
As the exploitation is relatively low, the convergence rate is usually low,
compared with that of traditional methods such as Newton-Raphson's methods.
 As a result, most metaheuristic algorithms can usually perform well for nonlinear problems,
including relatively tough optimization. However, the number of function evaluations
can be very high.

The role of crossover and mutation in exploration is rather subtle, while selection
as an exploitation mechanism can be simple and yet effective. However, it is still not
clear how the combination of crossover, mutation and selection can directly link
to the balance of exploration and exploitation. In fact, this is still an open question.
For example, in genetic algorithms, the probability of crossover can be as high as
0.95, while the mutation can be typically low in the range of 0.01 to 0.05.
In comparison with other algorithms, the exploration seems low, but genetic algorithms
have been proved to be very effective. On the other hand, mutation-related L\'evy flights
in cuckoo search can have a high exploration ability, and yet cuckoo search can converge
very quickly. At the same time, it is not clear what percentage of the search is
in exploration in the standard firefly algorithm, and yet it has been show that firefly
algorithm is very effective in dealing with multimodal, highly nonlinear problems.

Even in the standard particle swarm optimization, it is not clear what percentage of the
search iterations is in the exploration. The use of the current global best can be advantageous
and disadvantageous as well. The current global best may help to speed up the convergence,
but it may also lead to the wrong optimality if the current global best is obtained from
a biased set of samples drawn from a subspace where a local optimum (not the global optimum) is located.
All these suggest that it is still unknown how to achieve the optimal balance of exploration and exploitation
by tuning the combination of evolutionary operators.

In fact, a fine balance cannot be achieved by putting together all evolutionary operators
in a good way without tuning parameters. From experience, we know that the setting
or values of any algorithm-dependent parameters can affect the performance of an algorithm
significantly. In order to get good performance, we need to find the right values for
parameters. In other words, parameters need to be fine-tuned so that the algorithm can
perform to the best degree. Parameter tuning is still an active area of research \cite{Eiben}. An interesting development is a unified framework for self-tuning optimization algorithms \cite{Yang2013STA}.

Despite the importance of the above problems, little progress has been made.
On the contrary, there is some diversion in research efforts away from important problems.
Nature has evolved into millions of diverse species with
a diverse range of characteristics, but this does not mean that researchers should develop
millions of different algorithms, such as  the grass algorithm, leave algorithm, beatles algorithm,
sky algorithm,  universe algorithm, or even hooligan algorithm. Emphasis should focus on
solving important problems.

However, this does not mean new algorithms should not be developed at all. The research community
should encourage truly novel and efficient algorithms in terms of better evolutionary
operators and a better balance of exploration and exploitation.

\section{Conclusion}

Optimization algorithms based on swarm intelligence can have some distinct advantages
over traditional methods. By using theories of dynamical systems and self-organization
as well as the framework of Markov chains, we have provided a critical analysis of some
recent SI-based algorithms. The analysis has focus on the way of achieving exploration
and exploitation, and the basic components of evolutionary operators such as crossover,
mutation and selection of the fittest.

Through analysis, it has been found that most SI-based algorithms use mutation and selection
to achieve exploration and exploitation. Some algorithms use crossover as well, while most do not.
Mutation enables an algorithm to escape any local modes, while crossover provides good mixing to explore a
subspace more effectively, and thus more likely to lead to convergence. Selection provides a driving mechanism to
select the promising states or solutions.

The analysis also implies that there is room for improvement. Some algorithms such as PSO
may lack mixing and crossover, and therefore, hybridization may be useful to enhance its performance.

It is worth pointing out that the above analysis is based on the system behaviour for
continuous optimization problems, and it can be expected that these results are still valid
for combinatorial optimization problems. However, care should be taken for
combinatorial problems where neighbourhood may have different meaning, and, therefore,
the subspace concept may also be different. Further analysis and future studies may help to provide more
elaborate insight.

\end{document}